\documentclass[invmat, numbook, draft]{svjour}
\usepackage{times}
\usepackage{amsfonts}
\usepackage{amssymb,epic,eepic} %
\usepackage{english}
\author{ Igor Bjelakovi\'c \inst{1} \and Tyll Kr\"uger \inst{1,2}\and Rainer Siegmund-Schultze\inst{1}\and
  Arleta Szko\l a \inst{1}}
\institute {Technische Universit\"at Berlin
 \thanks{Fakult\"at II - Mathematik und Naturwissenschaften
 \\Institut f\"ur Mathematik MA 7-2
 \\Stra\ss e des 17. Juni 136
 \\10623 Berlin, Germany}, \email{\{igor, tkrueger, siegmund,
    szkola\}@math.tu-berlin.de} \and Universit\"at Bielefeld
 \thanks{Fakult\"at f\"ur Mathematik
 \\Universit\"atsstr. 25
 \\33619 Bielefeld, Germany}}
 \authorrunning{I. Bjelakovi\'c et al.}
\title{The Shannon-McMillan Theorem for Ergodic Quantum Lattice Systems}

\newcommand{\hr}{{\cal H}}

\newcommand{\nn}{{\mathbb N}}
\newcommand{\idn}{\mathbf{1}}
\newcommand{\zz}{{\mathbb Z}}
\newcommand{\eps}{{\varepsilon}}        

\newcommand{\x}{{\mathbf{x}}}
\newcommand{\n}{{\mathbf{n}}}
\newcommand{\kk}{{\mathbf{k}}}
\newcommand{\y}{{\mathbf{y}}}

\begin{document}
\maketitle

\begin{abstract}
We formulate and prove a quantum Shannon-McMillan theorem. The theorem demonstrates the
significance of the von Neumann entropy for translation invariant ergodic
quantum spin systems on $\zz^{\nu}$-lattices: the entropy gives the logarithm of
the essential number of eigenvectors of the system on large boxes. The
one-dimensional case covers quantum information sources and is basic
for coding theorems.
\end{abstract}

\section{Introduction}
The classical Shannon-McMillan theorem states that, for an ergodic shift system over a finite alphabet (spin chain, data stream, discrete 
stochastic process...) the 
distribution on a box of large volume $v$ is essentially carried by a 'typical set' of
approximate size $\exp (vh)$ with $h$ being the Kolmogorov-Sinai
dynamical entropy of the system, which in our case is identical to the
 mean (base $e$) Shannon entropy. Each element of this
typical set has a probability of order $\exp (-vh)$. The latter fact is usually
 called  \textsl{asymptotic equipartition property} (AEP). In order to
 simplify the arguments throughout this introduction we confine
 ourselves to the one-dimensional case. Thus instead of large boxes we
 consider $n$-blocks with $n$ large. In classical information
 theory the one-dimensional lattice systems model information sources.
 The AEP is of major importance for source coding (c.f. ~\cite{shields}). It implies the 
possibility of lossless transmission for an $n$-block message, using typically only
$nh / \log 2$ bits. The fact that on the typical set the probabilities of all elements are of
the same magnitude implies that better compression is (again typically) 
impossible.

Considering an ergodic quantum spin lattice system  (again, for simplicity, one-dimensional), we prove that the
$n$-block state is essentially carried by a typical subspace (of the entire Hilbert space) of dimension close
to $\exp (ns)$. Each one-dimensional
projector into that subspace has 
an expected value of order $\exp (-ns)$. Here $s$ is the mean von Neumann entropy of
the system. 
The reduction of this statement to the commutative setting is exactly the 'classical' 
theorem.
\\

Several attempts have been made to establish a quantum version of the 
Shannon-McMillan theorem for stationary quantum systems.
Due to the fact that in quantum  theory there exist several entropy notions (cf. \cite{benatti}), depending on e.g.
how the measuring process is incorporated, there is a number of
different approaches. The paper  \cite{king} of King and Le\'sniewski concerns the classical stochastic
process obtained from an ergodic quantum source by measurements of the
individual components and derives an asymptotic dimension of a
relevant subspace in terms of the classical mean entropy. This
relevant subspace however is in general not minimal, and
consequently not optimal for compression purposes. In particular, the AEP
is not fulfilled for this subspace. For the Bernoulli case (of
independent components) the King/Le\'sniewski result is optimal and consistent with
the earlier result by Jozsa and Schumacher \cite{jozsa}. 

In the paper \cite{petz} Hiai and Petz, using a considerably stronger
condition of complete ergodicity derived upper and lower bounds (in terms of
the von Neumann mean entropy) for the
asymptotic dimension
of a minimal relevant subspace. In view of our
results it turns out that the upper bound is sharp. However the results in \cite{petz} concern the
more general situation of estimating minimal projectors with respect
to a reference state in terms of relative entropy.
For the special cases of ergodic Gibbs quantum lattice
systems, \cite{gibbs}, respectively for ergodic algebraic states,
\cite{hiai}, Hiai and Petz could prove the
Shannon-McMillan theorem. For the Gibbs situation they even established the
convergence in the strong operator topology in the
corresponding GNS representation.

For further special classes of quantum systems, namely for a lattice
version of free bosons or fermions, Johnson and Suhov 
present a proof of a quantum Shannon-McMillan theorem, \cite{suhov}. The authors use
occupation numbers and their relation to the eigenvalues of Gibbs
ensemble density matrices to obtain a sequence of independent 
(not identically distributed) random variables. It turns out
that this classical process has a mean entropy which is precisely the
 mean von Neumann entropy of the corresponding quantum system.

Also in the context of quantum statistical mechanics Datta and Suhov \cite{datta} prove a
Shannon-McMillan theorem for some weakly non-stationary sequence of local quantum Gibbs states in
the neighborhood of a classical system. The authors use the quantum AEP to derive a reliable
compression/decompression scheme which compresses the quantum state according to the von Neumann entropy rate.
 
Neshveyev and St\o rmer considered the CNT-entropy in the context
of asymptotically abelian $C^{*}$-dynamical systems with locality in order to derive bounds for the asymptotic dimension of a typical
subspace, \cite{neshveyev}. They obtained a Shannon-McMillan result for the special case of tracial
states which covers the classical result. For our quantum lattice systems
however this essentially reduces to the situation of
Bernoulli states, where the individual components are even supposed to
be in the trace state.

Let us discuss, in an informal way, the main ideas behind the proof of the 
quantum Shannon-McMillan theorem presented here.

At the beginning, we follow 
mainly the approach of Hiai and Petz \cite{hiai}. It is well known that one can find a
classical subsystem (maximal abelian subalgebra of the entire
non-commutative algebra of observables) of the $n$-block
quantum system with Shannon entropy equal to the von Neumann entropy $S_{n}$
of the full $n$-block
quantum
state. Starting from this abelian $n$-block algebra we construct an
abelian quasilocal algebra. It is a subalgebra
of the entire quasilocal algebra of the infinite quantum lattice
system. Restricting the given quantum state to this classical infinite
lattice system we obtain a stochastic process. The $k$-th component of this
process is related to the $k$-th $n$-block of the
original system. This classical process is ergodic if one assumes for
instance \textsl{complete ergodicity} of the
original system as Hiai and Petz do in their paper. We know from classical Shannon-McMillan that the probability distribution of a large $k$-block of this classical
system (corresponding to a $kn$-block of the quantum system) has an essential
support of size $e^{k\cdot \widetilde{h}}$ (and entropy close to $%
k\cdot \widetilde{h}\leq k\cdot S_{n}$). Here $\widetilde{h}$ denotes
the mean entropy of the classical system, which is in general
smaller than $S_{n}$, due to possible correlations. For the quantum
context this translates into the existence of an essentially carrying subspace of dimension $%
e^{k\cdot \widetilde{h}}$ for the $kn$-block. With $s$ being the mean
von Neumann entropy, we arrive (for large $n$) at an estimate $%
kn\cdot s \; \sim \; k\cdot S_{n}\geq k\cdot \widetilde{h}\geq kn\cdot s$%
, and infer the existence of an (at most) $e^{k\cdot \widetilde{h}}\;\sim \;%
e^{kn\cdot s}$ -dimensional essentially carrying subspace ('typical
subspace') of a $kn$-block of the quantum system. If we would take into account
that the quantum system might have a still smaller ('better' from the point of view of resources needed)
typical subspace
(compared to what is delivered by the constructed classically restricted
system) we would have arrived at the result that $e^{n\cdot s}$ is an
\textsl{asymptotic upper bound }for the dimension of a typical subspace of
the quantum system. So far, we followed the ideas of Hiai and Petz  \cite{petz},
where under the assumption of complete ergodicity estimates for the typical
dimension were derived.

In fact, we can extend this result obtaining not only upper and lower bounds
but a limit assertion for the dimension of a typical subspace. Instead
of complete ergodicity we require only ergodicity,
thus establishing a quantum Shannon-McMillan theorem.

Beyond technicalities, this can be accomplished by

a) observing that we \textsl{do not }have to take into account the
possibility that the quantum system behaves 'better' than the classical one
(see some lines above): The very fact that the asymptotic von Neumann
entropy of a large $n$-block is $n\cdot s$ excludes the possibility of an
essentially carrying subspace of dimension significantly smaller than $%
e^{n\cdot s}$. In fact, the contribution to the entropy due to this subspace
would be smaller than $n\cdot s$. The non-typical part yields only a
 contribution of the order $o(n)$. This is due to a uniform integrability
argument (based upon Lemma \ref{class_resultat}) and

b) the condition of complete ergodicity can be avoided, since a quantum
analog of the following result concerning classical systems can be shown
(Theorem \ref{ergod_components}): if a dynamical system is ergodic with respect to the
transformation $T$, then, though it may be non-ergodic with respect to $T^{n}$ (for a
given $n\geq 2$), all of its $T^{n}$-ergodic components have the same
 entropy, and the number of ergodic components is at most $n$.

Note that the absence of complete ergodicity therefore reflects some intrinsic periodicity of the
  system. Periodicities naturally appear in a large variety
  of stationary (quantum) stochastic processes.

Finally let us mention that the techniques used here to derive the quantum Shannon-McMillan theorem
can be modified to prove a quantum version of Breiman's strengthening (in the sense of pointwise
convergence, \cite{breiman}) of the classical result. With the notion of an individual trajectory being
problematic in the quantum case, in the forthcoming paper \cite{bksss2} we give a reformulation of
Breiman's theorem, which has an immediate translation into the quantum situation.
%
\section{The Main Results}
Before we state our main theorems we shortly present the mathematical
description for the physical model of a $\nu$-dimensional quantum spin lattice
system, at the same time fixing notations.

The standard mathematical formalism is
introduced in detail e.g. in \cite{bratteli} and \cite{ruelle}. The $\nu$-dimensional infinitely extended lattice corresponds to the group
$\zz^{\nu}$. To each $\x \in \zz^{\nu}$ there is associated  an algebra
${\cal A}_{\mathbf{x}}$ of observables for a spin located at site $\x$. It is given by ${\cal
  A}_{\mathbf{x}}=\tau(\x){\cal A}$, where $\tau(\x)$ is an isomorphism and
${\cal A}$ is a finite dimensional $C^{*}$-algebra with identity. The local algebra ${\cal
  A}_{\Lambda}$ of observables for the finite subset $\Lambda \subset
\zz^{\nu}$ is given by ${\cal A}_{\Lambda}:= \bigotimes_{\x \in \Lambda} {\cal
  A}_{\x}$. The infinite lattice system is constructed from the finite subsets
$\Lambda \subset \zz^{\nu}$. It is described by the quasilocal $C^{*}$-algebra ${\cal
  A}^{\infty}$, which is defined as the operator norm closure of $\bigcup_{\Lambda \subset \zz^{\nu}}{\cal
  A}_{\Lambda}$.\\
A state of the
infinite spin system is given by a normed positive functional $\Psi$ on
${\cal A}^{\infty}$. It corresponds one-to-one to a consistent family of states
$\{\Psi^{(\Lambda)}\}_{\Lambda \subset \zz^{\nu}}$, where each
$\Psi^{(\Lambda)}$ is the restriction of $\Psi$ to the finite dimensional subalgebra ${\cal
  A}_{\Lambda}$ of ${\cal A}^{\infty}$ and consistency means that $\Psi^{(\Lambda)}=\Psi^{(\Lambda^{'})}\upharpoonright{\cal
  A}_{\Lambda}$ for $\Lambda \subset
\Lambda^{'}$. This one-to-one correspondence reflects the fact, that the state of the entire
spin lattice system is assumed to be determined by the expectation
values of all observables on finite subsystems $\Lambda$. Actually, it is
sufficient to consider boxes only. For each $\Psi^{(\Lambda)}$ there exists a unique density operator
$D_{\Lambda}\in {\cal A}_{\Lambda}$,
such that $\Psi^{(\Lambda)}(a)=\textrm{tr}_{\Lambda}D_{\Lambda}a,\ a
\in {\cal A}_{\Lambda}$ and $\textrm{tr}_{\Lambda}$ is the trace on
${\cal A}_{\Lambda}$. By ${\cal S}({\cal A}^{\infty})$ we denote the state space of ${\cal
  A}^{\infty}$.\\ Every $\mathbf{x} \in \zz^{\nu}$ defines a
translation of the lattice and induces an automorphism $T(\mathbf{x})$
on ${\cal
  A}^{\infty}$, which is a canonical extension of the isomorphisms for
finite $\Lambda \subset \zz^{\nu} $:
\begin{eqnarray}
T(\x): {\cal A}_{\Lambda} &\longrightarrow& {\cal A}_{\Lambda + \x}
\nonumber \\
a &\longmapsto& \left( \bigotimes_{\mathbf{z} \in
  \Lambda}T_{\mathbf{z}}(\x) \right) a,\nonumber
\end{eqnarray}
where $T_{\mathbf{z}}(\x):= \tau(\x)\tau^{-1}(\mathbf{z})$.
Then $\{T(\mathbf{x})\}_{\mathbf{x} \in
  \zz^{\nu}}$ is an action of the translation group $\zz^{\nu}$ by
automorphisms on ${\cal A}^{\infty}$.\\
Let  $G$ be any subgroup of $\zz^{\nu}$ and denote by ${\cal T}({\cal
  A}^{\infty}, G)$ the set of states, which are invariant under the
translations associated with $G$:
\begin{eqnarray}
{\cal T}({\cal A}^{\infty}, G):= \{ \Psi \in {\cal S}({\cal A}^{\infty})
| \Psi \circ T(\x) = \Psi, \ \forall \mathbf{x} \in G \}.  \nonumber
\end{eqnarray}
We will be concerned in particular with the space of $\zz^{\nu}$-invariant, i.e. stationary states. For simplicity we introduce the
abbreviation ${\cal T}({\cal
  A}^{\infty})={\cal T}({\cal
  A}^{\infty}, \zz^{\nu})$. Clearly, ${\cal T}({\cal
  A}^{\infty}) \subset{\cal T}({\cal
  A}^{\infty}, G)$ for any proper subgroup $G$ of $\zz^{\nu}$. \\ For
$\mathbf{n}=(n_{1}, \dots,n_{\nu})\in \nn^{\nu}$ we denote by $\Lambda
(\mathbf{n})$ the box in $\zz^{\nu}$ determined by
\begin{eqnarray}
\Lambda (\mathbf{n}):= \{(x_{1}, \dots, x_{\nu}) \in \zz^{\nu}| \ x_{i} \in
\{0, \dots, n_{i}-1\},\ i \in \{1, \dots,\nu \}\},   \nonumber
\end{eqnarray}
and for $n \in \nn$ the hypercube $\Lambda(n)$ is given by
\begin{eqnarray}
\Lambda (n):= \{ \x \in \zz^{\nu}| \ \x \in
\{0, \dots, n-1\}^{\nu}\}.
\end{eqnarray}
In the following we simplify notations by defining ${\cal A}^{(\mathbf{n})}:={\cal
  A}_{\Lambda(\mathbf{n})}$ and
$\Psi^{(\mathbf{n})}:=\Psi^{(\Lambda(\mathbf{n}))}$, for $\mathbf{n}
\in \nn^{\nu}$ (respectively for $n \in \nn$).\\
Obviously, because of the translation invariance any $\Psi \in {\cal
  T}({\cal A}^{\infty})$ is uniquely defined by the family $\{\Psi^{(\mathbf{n})}\}_{\mathbf{n} \in \nn^{\nu}}$. \\
The von Neumann entropy of $\Psi^{(\Lambda)}$ (cf. \cite{ohya}) is defined by:
\[S(\Psi^{(\Lambda)}):= -\textrm{tr}_{\Lambda}D_{\Lambda} \log
D_{\Lambda},\] where $\textrm{tr}_{\Lambda} $ denotes the trace on
${\cal A}_{\Lambda}$. For $\textrm{tr}_{\Lambda(\n)}$ we will write $\textrm{tr}_{\n}$. \\ It is well
known that for every $\Psi \in {\cal T}({\cal A}^{\infty})$ the limit 
\begin{eqnarray}
s(\Psi):=\lim_{\Lambda(\mathbf{n})
  \nearrow \nn^{\nu}} \frac{1}{\#(\Lambda(\mathbf{n}))}S(\Psi^{(\mathbf{n})}) \nonumber
\end{eqnarray}
exists. We call it the mean (von Neumann) entropy.
Let $l \in \nn$ and consider the subgroup $G_{l}:= l \cdot
\zz^{\nu}$. For a $G_{l}$-invariant state $\Psi$ we define the mean entropy
with respect to $G_{l}$ by
\begin{eqnarray}\label{G_l_mean_entropy}
s(\Psi, G_{l}):=\lim_{\Lambda(\mathbf{n}) \nearrow \nn^{\nu}} \frac{1}{\#(\Lambda(\mathbf{n}))}S(\Psi^{(l\cdot\mathbf{n})}).
\end{eqnarray}
Observe that if the state $\Psi$ is $\zz^{\nu}$-invariant then we have
the relation $s(\Psi, G_{l})= l^{\nu}\cdot s(\Psi)$. 
Further note that in the commutative case the mean entropy $s$ coincides with
the Kolmogorov-Sinai dynamical entropy $h$. \\
${\cal T}({\cal A}^{\infty}, G)$ is a convex, (weak-*) compact subset of ${\cal
  S}({\cal A}^{\infty})$ for any subgroup $G$ of $\zz^{\nu}$. We denote by $\partial_{\textrm{ex}}{\cal T}({\cal
  A}^{\infty}, G)$ the set of extreme points of ${\cal
  T}({\cal A}^{\infty}, G)$. We refer to elements of $\partial_{\textrm{ex}}{\cal T}({\cal
  A}^{\infty}, G)$ as  $G$-ergodic states.  The elements of $\partial_{\textrm{ex}}{\cal T}({\cal
  A}^{\infty})$ are called ergodic states.

Now we formulate our main result as a generalization of
the Shannon-McMillan theorem to the case of non-commutative quasilocal
$C^{*}$-algebras ${\cal A}^{\infty}$ constructed from ${\cal A}$, a finite dimensional $C^{*}$-algebra with identity.
\begin{theorem}[Quantum Shannon-McMillan Theorem]\label{SM}
Let \(\Psi\) be an ergodic state on \( {\cal A}^{\infty}\) with
 mean entropy \(s(\Psi)\). Then
for all \( \delta > 0\) there is an \( \mathbf{N}_{ \delta} \in \nn^{\nu}\)
such that for all $\n \in \nn^{\nu}$ with \(\Lambda(\n) \supseteq \Lambda(\mathbf{N}_{ \delta})\) there exists an
orthogonal projector \(p_{\n}(\delta) \in {\cal A}^{(\n)}\) such that
\begin{enumerate}
\item
\(\Psi^{(\n)}(p_{\n}(\delta))\geq 1-\delta\),
\item
for all minimal projectors \(0\neq p \in {\cal A}^{(\n)}\) with \(p \leq p_{\n}(\delta)\)
\[e^{-\#(\Lambda(\n))(s(\Psi)+ \delta)} < \Psi^{(\n)}(p) <
e^{-\#(\Lambda(\n))(s(\Psi)- \delta)},\]
\item
\( e^{\#(\Lambda(\n))(s(\Psi)- \delta)} < \textup{tr}_{\n}(p_{\n}(\delta)) < e^{\#(\Lambda(\n))(s(\Psi)+ \delta)}.\)
\end{enumerate}
\end{theorem}
\textbf{Remark:}
\emph{In complete analogy to the classical case the above theorem, in particular
item 2., expresses the AEP in a
quantum context. As will be seen in the proof of this theorem the typical subspace given
by the projector $p_{\n}(\delta)$ can be chosen as the linear hull of the eigenvectors of
$\Psi^{(\n)}$ which have eigenvalues of order
$e^{-\#(\Lambda(\n))s(\Psi)}$.}\\ \\ 
In quantum information theory spin systems on one-dimensional
lattices model quantum information sources. In such a setting essentially carrying
subspaces of  minimal dimension are particularly important for
formulating and proving coding or compression theorems. 
As subspaces
of probability close to 1 they are the relevant subspaces in the sense
that the expectation values of any observables restricted to these subspaces are
almost equal to the corresponding ones on the entire space. On the other hand the
minimal dimension allows an economical use of resources (qubits) needed for
quantum information storage and transmission.
Define for $\eps \in (0,1)$ and $\mathbf{n} \in \nn^{\nu} $
\begin{eqnarray}
\beta_{\eps, \mathbf{n}}(\Psi):= \min \{ \log ( \textrm{tr}_{\mathbf{n}}\ q) | \ q \in {\cal
  A}^{(\mathbf{n})} \ \textrm{projector}, \ \Psi^{(\mathbf{n})}(q) \geq 1- \eps\}, \nonumber
\end{eqnarray}
Our next statement, which appears already as a conjecture in the
papers \cite{hiai} and \cite{petz} of Petz and Hiai,
is strongly related to the previous
one. In fact, it can be considered as a
reformulation especially for possible applications to (quantum) data compression.
\begin{proposition}\label{beta=s}
Let \(\Psi\) be an ergodic state on \( {\cal A}^{\infty}\) with
 mean entropy \(s(\Psi)\). Then for every $\eps
\in (0,1)$
\begin{eqnarray}\label{lim_beta=s}
\lim_{\Lambda(\n) \nearrow \nn^{\nu}}\frac{1}{\#(\Lambda(\n))} \beta_{\eps, \n}(\Psi)=s(\Psi).
\end{eqnarray}
\end{proposition}
In the proof of the quantum Shannon-McMillan theorem this proposition
is an intermediate result.
\section{Proofs of the Main Results}
A basic tool for the proof of the Shannon-McMillan theorem under the general assumption of ergodicity is the structural
assertion Theorem \ref{ergod_components}. It is used to circumvent the complete ergodicity assumption made by Hiai and Petz.
The substitution of the quantum system by a classical approximation on large boxes leads to an ergodicity problem for these classical approximations.
Theorem \ref{ergod_components} combined with the subsequent lemma allow to control not only the mean (per site limit)
entropies of the ergodic components (with respect to the subshift generated by a large box), but also to cope with the obstacle
that some of these components
might have an atypical entropy on this large but finite box. Using these prerequisites, we prove Lemma
\ref{lim_beta<s} which is the extension of the Hiai/Petz upper bound result to ergodic states. Finally, from the simple
probabilistic argument expressed in Lemma \ref{class_resultat} we infer that the upper bound is really a limit.
\begin{theorem}\label{ergod_components}
Let $\Psi \in \partial_{\textrm{ex}} {\cal T} ({\cal A}^{\infty})$. Then for
every subgroup $G_{l}:=l \cdot \zz^{\nu}$, with $l>1$ an integer, there exists a ${\kk}(l) \in {\nn}^{\nu}$ and a unique
convex  decomposition of $\Psi $ into $G_{l}$-ergodic
 states $ \Psi_{\x}$:

\begin{eqnarray}\label{ergodic_decomposition}
\Psi =\frac{1}{\#(\Lambda ({\kk}(l)))}\sum_{\x\in \Lambda ({\kk}(l))} \Psi_{\x}.
\end{eqnarray}
The $G_{l}$-ergodic decomposition (\ref{ergodic_decomposition}) has the following properties:
\begin{enumerate}
\item $k_{j}(l) \leq l $ and $k_{j}(l)|l $ for all \(j\in
\{1,\dots ,\nu\}\)

\item \(\{\Psi_{\x}\}_{\x\in \Lambda ({\kk}(l))}= \{\Psi_{0} \circ T(-\x)\}_{\x\in \Lambda
  ({\kk}(l))} \)

\item For every \(G_{l}\)-ergodic state \(\Psi_{\x}\) in the
convex decomposition (\ref{ergodic_decomposition}) of \(\Psi \) the mean entropy  with respect to \(G_{l}\), \(s(\Psi_{\x}, G_{l})\), is
equal to the mean entropy
\(s(\Psi, G_{l})\), i.e. 
\begin{eqnarray}\label{equal_entropies}
s(\Psi_{\x}, G_{l}) = s(\Psi, G_{l}) 
\end{eqnarray}
for all $\x\in \Lambda ({\kk}(l))$.
\end{enumerate}
\end{theorem}
\textbf{Proof of Theorem \ref{ergod_components}:}
Let $(\hr_{\Psi},\pi_{\Psi}, \Omega_{\Psi}, U_{\Psi})$ be the GNS
representation of the $C^{*}$-dynamical system $({\cal A}^{\infty}, \Psi,
{\zz}^{\nu})$. $U_{\Psi}$ is the unitary representation of ${\zz}^{\nu}$ on $\hr_{\Psi}$. It
satisfies for every $\x\in{\zz}^{\nu}$:
\begin{eqnarray}
U_{\Psi}(\x)\Omega_{\Psi}&=&\Omega_{\Psi}, \label{omega_u_inv} \\
U_{\Psi}(\x)\pi_{\Psi}(a)U_{\Psi}^{*}(\x)&=&\pi_{\Psi}(T(\x)a), \qquad
\forall a \in {\cal A}^{\infty}.  \label{Ta}
\end{eqnarray}
Define 
\begin{eqnarray}
{\cal N}_{\Psi, G_{l}}&:=& \pi_{\Psi}({\cal A}^{\infty}) \cup
 U_{\Psi}(G_{l}), \\ {\cal P}_{\Psi, G_{l}}&:=& \{P \in
{\cal N}'_{\Psi, G_{l}} |\ P= P^{*}=P^{2}\}.
\end{eqnarray} 
By $'$ we denote the commutant. Observe that ${\cal N}_{\Psi, G_{l}}$  is selfadjoint. Thus  ${\cal
  N}'_{\Psi, G_{l}}$ (as the commutant of a selfadjoint set) is a
von Neumann algebra. Further it is a known result that ${\cal
  N}'_{\Psi, G_{l}}$ is abelian, (cf. Proposition 4.3.7. in
\cite{bratteli} or Lemma IV.3.4 in ~\cite{israel}). \\ \\ Consider some $l>1$ such that  $\Psi \notin \partial_{\textrm{ex}} {\cal T}({\cal
  A}^{\infty}, G_{l})$. (If there is no such $l$ the statement of the
  theorem is trivial.) Then
\begin{eqnarray}\label{pkt_1}
{\cal P}_{\Psi, G_{l}} \backslash \{0, \idn\} \not= \emptyset.
\end{eqnarray}
In fact, $\Psi \notin \partial_{\textrm{ex}} {\cal T}({\cal A}^{\infty}, G_{l})$ is
equivalent to the reducibility of ${\cal N}_{\Psi, G_{l}}$,
(cf. Theorem 4.3.17 in \cite{bratteli}). This means that there is a
non-trivial closed subspace of $\hr_{\Psi}$ invariant under the action
of $\pi_{\Psi}(\mathcal{A}^{\infty})$ and $U_{\Psi}(G_{l})$. Let $P$ be 
the projection on this subspace and $
P^{\perp}=\idn -P$ . Then $P,P^{\perp} \notin \{ 0,  \idn\}$ and of course $P$
and $ P^{\perp}$ are contained in ${\cal N}'_{\Psi, G_{l}}$. Thus
(\ref{pkt_1}) is clear. \\ \\ Let $I$ be a countable
index set. An implication of the ${\zz}^{\nu}$-ergodicity of the $G_{l}$-invariant $\Psi$
is the following:
\begin{eqnarray} \label{pkt_2}
\{Q_{i} \}_{i \in I}\ \textrm{orthogonal}\ \idn - \textrm{decomposition in}\
{\cal N}'_{\Psi, G_{l}} \Longrightarrow |I| \leq l^{\nu}.
\end{eqnarray}
To see (\ref{pkt_2}) observe at first that  for any \(Q \in  {\cal P}_{\Psi,
  G_{l}}\backslash \{ 0 \}\) the projector $U_{\Psi}(\x) Q
U^{*}_{\Psi}(\x) $, $\x\in \Lambda (l)$, belongs to the abelian algebra \({\cal
  N}'_{\Psi, G_{l}}\), namely
\begin{eqnarray}
\pi_{\Psi}(a)U_{\Psi}(\x) Q
U^{*}_{\Psi}(\x) &=& U_{\Psi}(\x)\pi_{\Psi}(T(-\x)a) Q U^{*}_{\Psi}(\x) \quad
(\textrm{by } (\ref{Ta}))
\nonumber \\ &=& U_{\Psi}(\x)Q \pi_{\Psi}(T(-\x)a) U^{*}_{\Psi}(\x) \quad
\nonumber \\ &=& U_{\Psi}(\x)Q U^{*}_{\Psi}(\x)\pi_{\Psi}(a)  \quad
(\textrm{by } (\ref{Ta}))\nonumber
\end{eqnarray}
holds for every \(a\in{\cal A}^{\infty}\) and \([U_{\Psi}(\mathbf{y}),U_{\Psi}(\x) Q
U^{*}_{\Psi}(\x)]=0\) is 
obvious by  \([U_{\Psi}(\mathbf{y}),U_{\Psi}(\x)]=0\) for all $\mathbf{y}\in G_{l}$ and $\x\in {\zz}^{\nu}$. Thus $\{U_{\Psi}(\x) Q
U^{*}_{\Psi}(\x)\}_{\x\in \Lambda (l)}$ is a family of mutually
commuting projectors. The Gelfand isomorphism represents the projectors
$U_{\Psi}(\x) Q U^{*}_{\Psi}(\x)$ as
 continuous characteristic functions $1_{Q_{x}}$ on some compact (totally
 disconnected) space. Define
\(\bar{Q} := \bigvee
_{\x\in\Lambda (l)} U_{\Psi}(\x) Q U^{*}_{\Psi}(\x)\), which has the
representation as $\bigvee_{\x \in \Lambda (l)}
1_{Q_{\x}}=1_{\bigcup_{\x\in \Lambda (l)}Q_{\x}}$. Note that \(\bar{Q}\)
is 
invariant under the action of $U_{\Psi}({\zz}^{\nu})$. From the 
${\zz}^{\nu}$-ergodicity of $\Psi$ we deduce that
$\bar{Q} = \idn$. If
 we translate back the finite subadditivity of probability measures to
the expectation values of the projectors $U_{\Psi}(\x) Q U^{*}_{\Psi}(\x)$
we obtain: 
\begin{eqnarray}
1= \langle \Omega_{\Psi}, \bar{Q} \Omega_{\Psi} \rangle &\leq&
\sum_{\x\in \Lambda (l)}   \langle \Omega_{\Psi},U_{\Psi}(\x) Q
U^{*}_{\Psi}(\x)\Omega_{\Psi} \rangle \nonumber \\ &=& l^{\nu} \cdot \langle \Omega_{\Psi},
Q \Omega_{\Psi} \rangle \quad (\textrm{by } (\ref{omega_u_inv})).
\nonumber   
\end{eqnarray}

Thus (\ref{pkt_2}) is clear.\\
Combining the results (\ref{pkt_1}) and (\ref{pkt_2}) we get the
existence of an orthogonal $\idn$-decomposition $\{P_{i}\}_{i=0}^{n_{l}-1}$ in $n_{l}
\leq l^{\nu}$ minimal projectors $ P_{i} \in {\cal P}_{\Psi,
  G_{l}} \backslash \{0,\idn \}$. Here we use the standard definition of minimality:
\begin{eqnarray}
P \ \textrm{minimal projector in} \ {\cal N}'_{\Psi, G_{l}}\ 
&:\Longleftrightarrow& 0 \not= P \in {\cal P}_{\Psi, G_{l}} \ \textrm{and}
\nonumber \\ &  Q \leq P &\ \Rightarrow Q =P \quad \forall Q
\in {\cal P}_{\Psi, G_{l}}\backslash \{0\} \nonumber
\end{eqnarray}
The abelianness of ${\cal N}'_{\Psi, G_{l}}$ implies the uniqueness of the orthogonal
$\idn $-decomposition $\{P_{i}\}_{i=0}^{n_{l}-1}$. Further it follows that $\{P_{i}\}_{i=0}^{n_{l}-1}$ is a generating subset for
${\cal P}_{\Psi, G_{l}}$ in the following sense:
\begin{eqnarray}\label{generating_set}
Q \in {\cal
  P}_{\Psi, G_{l}}\ \Longrightarrow\ \exists \{P_{i_{j}}\}_{j =0}^{s \leq n_{l}-1}
\subset {\cal
  P}_{\Psi, G_{l}}\ \textrm{ such that } Q = \sum_{j=0}^{s}P_{i_{j}}.
\end{eqnarray}
Define $p_{i}:=
\langle \Omega_{\Psi}, P_{i} \Omega_{\Psi} \rangle$ and order the
minimal projectors $P_{i}$ such that
\begin{eqnarray}\label{ordering}
p_{0} \leq p_{i}, \qquad \forall i \in \{1, \dots , n_{l}-1\}.
\end{eqnarray}
Let
\[
G(P_{0}):=\{\x\in {\zz}^{\nu}|\ U(\x)P_{0}U^{*}(\x)=P_{0}\}.
\]
Note that $G(P_{0})$ is a subgroup of ${\zz}^{\nu}$ and contains
$G_{l}$, since $P_{0} \in {\cal P}_{\Psi, G_{l}}$.  This leads to the representation
\[ G(P_{0})=\bigoplus_{j=1}^{\nu}k_{j}(l)\zz, \quad \textrm {with }k_{j}(l)|l\quad \textrm{
  for all } j\in \{1,\dots ,\nu\},\]
where the integers $k_{j}(l)$ are given by
\[ k_{j}(l):=\min\{x_{j}|\ x_{j} \textrm{ is the j-th component of }\x\in
G(P_{0})\textrm{ and } x_{j}>0\}.\]
For $P_{0}$, as an element of ${\cal P}_{\Psi, G_{l}}$,
$\{U_{\Psi}(\x) P_{0} U^{*}_{\Psi}(\x)\}_{\x\in \Lambda ({\kk}(l))} \subseteq {\cal
  P}_{\Psi, G_{l}}$ for ${\kk}(l)=(k_{1}(l),\dots ,k_{\nu}(l))$.
 Thus by (\ref{generating_set}) each $U_{\Psi}(\x)
P_{0} U^{*}_{\Psi}(\x) $, $\x\in \Lambda ({\kk}(l))$, can be represented as a sum of minimal
projectors. But then by linearity of the expectation values and the
assumed ordering (\ref{ordering}) each $U_{\Psi}(\x)
P_{0} U^{*}_{\Psi}(\x) $ must be a minimal projector for $\x\in \Lambda
({\kk}(l))$.
Otherwise there
would be a contradiction to $\langle \Omega_{\Psi},U_{\Psi}(\x) P_{0}
U^{*}_{\Psi}(\x) \Omega_{\Psi} \rangle = p_{0}$. Consequently $\{U_{\Psi}(\x) P_{0}
U^{*}_{\Psi}(\x)\}_{\x\in \Lambda ({\kk}(l))} \subseteq \{P_{i}\}_{i=0}^{n_{l}-1}$. 
Consider $\bar P_{0} =
\sum_{\x\in \Lambda ({\kk}(l))} U_{\Psi}(\x) P_{0} U^{*}_{\Psi}(\x) $. Obviously $\bar
P_{0}$ is invariant under the action of $U_{\Psi}({\zz}^{\nu})$ and because of the
${\zz}^{\nu}$-ergodicity of $\Psi$
\[ \bar P_{0} = \idn.\]
It follows by the uniqueness of the orthogonal $\idn$-decomposition
\[ \{U_{\Psi}(\x)
P_{0}U^{*}_{\Psi}(\x)\}_{\x\in \Lambda ({\kk}(l))}=\{P_{j}\}_{j=0}^{n_{l}-1}. \]
Obviously $n_{l}=\#(\Lambda(\kk(l)))$ and for each $P_{i}$,
$i\in \{0, \dots ,n_{l}-1\}$, there is only one $\x\in \Lambda ({\kk(l)})$ such that
\begin{eqnarray}\label{upu}
P_{i}=U_{\Psi}(\x) P_{0}U^{*}_{\Psi}(\x)=:P_{\x}.
\end{eqnarray}
It follows
\(
p_{i}=p_{0}\) for all $i \in \{0, \dots, n_{l}-1 \}$ and hence 
\begin{eqnarray}
p_{i}=\frac{1}{n_{l}}=\frac{1}{\#(\Lambda(\kk(l)))}, \qquad i \in \{0, \dots, n_{l}-1\}. \nonumber
\end{eqnarray}
Finally, set for every $\x \in \Lambda ({\kk}(l))$
\begin{eqnarray}
\Psi_{\x}(a) &:=& \#(\Lambda ({\kk}(l))) \langle
\Omega_{\Psi}, P_{\x} \pi_{\Psi}(a) \Omega_{\Psi} \rangle, \qquad a \in {\cal A}^{\infty}. \nonumber 
\end{eqnarray}
From (\ref{upu}), (\ref{omega_u_inv}) and (\ref{Ta}) we get
\begin{eqnarray}
\Psi_{\x}(a) &=& \#(\Lambda ({\kk}(l))) \langle
\Omega_{\Psi}, P_{\x} \pi_{\Psi}(a) \Omega_{\Psi} \rangle  \nonumber \\
 &=&  \#(\Lambda ({\kk}(l))) \langle
\Omega_{\Psi}, P_{0} \pi_{\Psi}(T(-\x)a) \Omega_{\Psi} \rangle
\nonumber \\ &=& \Psi_{0}(T(-\x)a), \qquad a \in {\cal A}^{\infty}, \nonumber
\end{eqnarray}
hence
\begin{eqnarray*}
\frac{1}{\#(\Lambda(\kk(l)))}\sum_{\x \in \Lambda ({\kk}(l))} \langle
\Omega_{\Psi}, P_{\x} \pi_{\Psi}(a) \Omega_{\Psi} \rangle 
&=&  \langle
\Omega_{\Psi}, (\sum_{\x\in\Lambda ({\kk}(l))} P_{\x}) \pi_{\Psi}(a) \Omega_{\Psi}
\rangle \\ &=& \Psi (a).
\end{eqnarray*}
 Thus we arrive at the  convex
 decomposition of $\Psi$:
\begin{eqnarray}
\Psi=\frac{1}{\#(\Lambda ({\kk}(l)))}\sum_{\x\in \Lambda ({\kk}(l))}
\Psi_{0}\circ T(-\x) . \nonumber
\end{eqnarray}
By construction this is a $G_{l}$-ergodic decomposition of
\(\Psi\).
%
%
%
%
It remains to prove the fact that the
mean entropies with
respect to the lattice $G_{l}$ are the same for all $G_{l}$-ergodic
components ${\Psi}_{\x}$.    

\textsl{Proof of item 3.}: It is a well known result that the quantum mean
entropy with respect to a given lattice $G_{l}$ is affine on the convex set of 
$G_{l}$-invariant states, (cf. prop. 7.2.3 in \cite{ruelle}). Thus to prove (\ref{equal_entropies}) it is sufficient to show:
\begin{eqnarray}
s(\Psi _{\mathbf{x}},G_{l})=s(\Psi _{0},G_{l}),\qquad \forall \mathbf{x}\in
\Lambda (\mathbf{k}(l)).\nonumber
\end{eqnarray}
By the definition of the mean entropy this is equivalent to the statement 
\begin{eqnarray}\label{abschaetzung}
|S(\Psi _{\mathbf{x}}^{(l\mathbf{n})})-S(\Psi _{0}^{(l\mathbf{n})})|=o(|%
\mathbf{n|)} \qquad \textrm{as }\mathbf{n\rightarrow \infty .}
\end{eqnarray}
This can be seen as follows: In view of the definition of \ $\Psi _{\mathbf{x}}^{(l\mathbf{n})}$ we
have\  $S(\Psi _{\mathbf{x}}^{(l\mathbf{n})})=S(\Psi _{\mathbf{x}}^{(\Lambda
(l\mathbf{n}))})=S(\Psi _{0}^{(\Lambda
(l\mathbf{n})-\mathbf{x})})$. We introduce the box $\tilde{\Lambda}$
being concentric with $\Lambda (l\mathbf{n),}$ with all edges enlarged by $l$
on both directions, i.e. an $l$-neighborhood of  $\Lambda (l\mathbf{n).}$
The two expressions $S(\Psi _{\mathbf{x}}^{(l\mathbf{n})})$ and $S(\Psi
_{0}^{(l\mathbf{n})})$ are von Neumann entropies of the restrictions
of $\Psi _{0}^{(\widetilde{\Lambda})}$ to the smaller sets $\Lambda (l\mathbf{n})$ and $%
\Lambda (l\mathbf{n})-\mathbf{x}$, respectively. 
On the other hand we consider the box $\widehat{\Lambda}$ being concentric with $\Lambda(l\n)$ with all edges shortened by $l$
at both sides. $S({\Psi}_{0}^{(\widehat{\Lambda})})$ is the von Neumann entropy of
${\Psi}_{0}^{(\Lambda(l\n))}$ and ${\Psi}_{\x}^{(\Lambda(l\n))}$
after their restriction to the set $\widehat{\Lambda}$.
$ S({\Psi}_{\x}^{(\Lambda(l\n))})$ and $S({\Psi}_{0}^{(\Lambda(l\n))})$ 
can be estimated
simultaneously using the subadditivity of the von Neumann
entropy 
\begin{eqnarray}
S(\Psi _{0}^{(\widetilde{\Lambda })})-\log \textrm{tr}_{\tilde{\Lambda}\backslash \Lambda (l%
\mathbf{n})}\mathbf{1} \leq S(\Psi^{(l \mathbf{n})}_{\natural}) \leq S(\Psi _{0}^{%
(\widehat{\Lambda })})+ \log \textrm{tr}_{\Lambda(l\n)\backslash \widehat{\Lambda}}\mathbf{1}, \nonumber
\end{eqnarray}
where $\natural\in\{ \mathbf{x}, 0\}$.
Thus (\ref{abschaetzung}) is
immediate.$\qquad \qed$ \\ 

In order to simplify our notation in the next lemma we introduce some 
abbreviations. We choose a
positive integer $ l $ and consider the decomposition of $\Psi \in
{\partial}_{ ex } \mathcal{T} (\mathcal{ A^{ \infty } })$  
into states ${\Psi}_{\x}$ being ergodic with respect to the
action of $G_{l}$, i.e. $\Psi=\frac{1}{\#(\Lambda (\kk (l)))}\sum_{\x
  \in \Lambda (\kk (l))}{\Psi}_{\x}$. 
Then we set 
\[ s:=s(\Psi ,{\zz}^{\nu})=s(\Psi),\]
i.e. the mean entropy of the state $\Psi$ computed with respect to ${\zz}^{\nu}$. Moreover
we set 
\[s_{\x}^{(l)}:=\frac{1}{\#(\Lambda (l))}S({\Psi}_{\x}^{(\Lambda(l))})\quad
\textrm{ and }\quad s^{(l)}:=\frac{1}{\#(\Lambda
  (l))}S({\Psi}^{(\Lambda (l))}).\]
From the previous lemma we know that
\begin{equation} \label{eq:el}
s(\Psi_{\x} ,G_{l})=s(\Psi ,G_{l})=l^{\nu} \cdot s(\Psi),\qquad \forall
\x\in \Lambda (\kk (l)).
\end{equation}
For $\eta > 0$ let us introduce the following set
\begin{equation}\label{eq:a}
 A_{l,\eta}:=\{\x \in \Lambda (\kk (l))| \
  s_{\x}^{(l)}\ge s+\eta\}.
\end{equation}
By \(A_{l,\eta}^{c}\) we denote its complement.
The following lemma states that the density of $G_{l}$-ergodic
components of $\Psi$ which have too large entropy on the box of side length
$l$ vanishes asymptotically in $l$.
\begin{lemma}\label{asimpt}
If $\Psi$ is a ${\zz}^{\nu}$-ergodic state on $ {\mathcal{A}}^{\infty}$, then 
\[ \lim_{l \rightarrow \infty}\frac{\# A_{l,\eta}}{\# \Lambda (\kk (l))}=0\]
holds for every $\eta >0$.
\end{lemma}
\textbf{Proof of Lemma \ref{asimpt}:} We suppose on the contrary that there is some $\eta_{0} >0$
such that $\limsup_{l}\frac{\# A_{l,\eta_{0}}}{\# \Lambda (\kk
  (l))}=a>0$. 
Then there exists a subsequence $ (l_{j})$ with the property  
\[ \lim_{j \rightarrow \infty}\frac{\#
  A_{l_{j},\eta_{0}}}{\# \Lambda (\kk (l_{j}))}=a.\]
By the concavity of the von Neumann entropy we obtain
\begin{eqnarray} \#\Lambda (\kk (l_{j})) \cdot s^{(l_{j})} & \geq &
\sum_{\x \in \Lambda (\kk
  (l_{j}))}s_{\x}^{(l_{j})}\nonumber \\
&=&\sum_{\x\in
 A_{l_{j},\eta_{0}}}s_{\x}^{(l_{j})}+\sum_{\x\in
 A_{l_{j},\eta_{0}}^{c}}s_{\x}^{(l_{j})} \nonumber \\
& \geq & \# A_{l_{j},\eta_{0}} \cdot (s+\eta_{0})+ \#
  A_{l_{j},\eta_{0}}^{c} \cdot \min_{\x\in
  A_{l_{j},\eta_{0}}^{c}}s_{\x}^{(l_{j})}.
\nonumber
\end{eqnarray}
Here we made use of (\ref{eq:a}) at the last step. Using that for the
mean entropy holds
\[s({\Psi}_{\x} ,G_{l})=\lim_{\Lambda(\mathbf{m})\nearrow \nn^{\nu}}\frac{1}{\#\Lambda
    (\mathbf{m})}S({\Psi}_{\x}^{(l\mathbf{m})})=\inf_{\Lambda
    (\mathbf{m})}
\frac{1}{ \# \Lambda (\mathbf{m})}S({\Psi}_{\x}^{(l\mathbf{m})})\] 
we obtain a further estimation for the second term on the right hand
side:
\begin{eqnarray*}
\#
  A_{l_{j},\eta_{0}}^{c} \cdot \min_{\x\in
  A_{l_{j},\eta_{0}}^{c}}s_{\x}^{(l_{j})}
&\geq & \#
  A_{l_{j},\eta_{0}}^{c} \cdot \min_{\x\in
  A_{l_{j},\eta_{0}}^{c}}\frac{1}{l_{j}^{\nu}}s({\Psi}_{\x}
,G_{l_{j}})
\\
& = & \#
  A_{l_{j},\eta_{0}}^{c} \cdot s(\Psi) \qquad \textrm{
  (by(\ref{eq:el}))}.
\end{eqnarray*}
 After dividing $\#\Lambda (\kk (l_{j})) \cdot s^{(l_{j})} \geq \#
A_{l_{j},\eta_{0}} \cdot (s+\eta_{0})+ \#
  A_{l_{j},\eta_{0}}^{c} \cdot s(\Psi)$  by $\#\Lambda (\kk (l_{j}))$
  and taking limits we arrive at the following contradictory inequality:
\[ s\geq a(s+\eta_{0})+(1-a)s=s+a\eta_{0}>s.\]
So, $a=0$.$\qquad \qed$ \newline
\begin{lemma}\label{lim_beta<s}
Let $\Psi$ be an ergodic state on $\mathcal{ A^{ \infty } }$. Then
for every $\varepsilon \in (0,1)$
\[ \limsup_{\Lambda(\n) \nearrow \nn^{\nu}}\frac{1}{\#\Lambda(\n)} \beta_{\varepsilon,\n}(\Psi)\leq s(\Psi).\]
\end{lemma}
\textbf{Proof of Lemma \ref{lim_beta<s}:}
 We fix $\varepsilon >0$ and choose arbitrary $\eta$,
$\delta>0$. 
Consider the $G_{l}$-ergodic decomposition 
 \[\Psi=\frac{1}{\#\Lambda (\kk(l))}\sum_{\x \in \Lambda (\kk(l))}\Psi_{x}\]
of $\Psi$ for integers $l\geq1$. By Lemma \ref{asimpt} there is 
an integer $L\geq 1$ such that for any $l\geq L$ 
\[\frac{\varepsilon}{2} \geq \frac{1}{\# \Lambda (\kk(l))}\#{A_{l,\eta}}\geq 0\]
holds, where $A_{l,\eta}$ is defined by (\ref{eq:a}). This inequality implies 
\begin{equation} \label{eq:typical}
 \frac{1}{\#\Lambda (\kk(l))}\#{A_{l,\eta}^{c}} \cdot (1-\frac{\varepsilon}{2})\geq 1-\varepsilon. 
\end{equation}
On the other hand by
\[S(\Psi^{(\n)})=\inf\{S(\Psi^{(\n)}\!\upharpoonright\mathcal{B})|\ \mathcal{B}\textrm{
  maximal abelian }C^{\ast}-\textrm{subalgebra of }\mathcal{A}^{(\n)}\}\]
(cf. Theorem 11.9 in~\cite{niel} and use the one-to-one correspondence between
 maximal abelian $\ast$-subalgebras and orthogonal
 $\mathbf{1}$-decompositions into minimal projectors contained in $\mathcal{A}^{(\n)}$) 
there exist maximal abelian $C^{\ast}$-subalgebras
$\mathcal{B}_{\x}$ of $ \mathcal{A}_{\Lambda (l)}$ with the property
\begin{equation}\label{eq:goodentr}
\frac{1}{\#\Lambda(l)}S(\Psi_{\x}^{(\Lambda(l))}\!\upharpoonright\mathcal{B}_{\x})<s(\Psi)+\eta, \qquad
\forall \ \x \in A_{l,\eta}^{c}.  
\end{equation}
We fix an $l \geq L$ and consider the abelian
quasi-local $C^{\ast}$-algebras $\mathcal{B}_{\x}^{\infty}$,
constructed with $\mathcal{B}_{\x}$,
as $C^{\ast}$-subalgebras of $\mathcal{A}^{\infty}$ and set
\[\ m_{\x}:=\Psi_{\x}\!\upharpoonright\mathcal{B}_{\x}^{\infty} 
\textrm{ and } m_{\x}^{(\n)}:=\Psi_{\x}\!\upharpoonright\mathcal{B}_{\x}^{(\n)}\]
for $\x\in A_{l,\eta}^{c}$ and $\n\in {\nn}^{\nu}$. The states $m_{\x}$ are
$G_{l}$-ergodic since they are restrictions of $G_{l}$-ergodic states
$\Psi_{\x}$ on a quasi-local algebra. This easily follows from Theorem
4.3.17. in \cite{bratteli}.
 Moreover, by the Gelfand isomorphism and Riesz representation
 theorem, we can identify the states $m_{\x}$ with
probability measures on corresponding (compact) maximal ideal spaces
 of $\mathcal{B}_{\x}^{\infty}$.
By commutativity and finite dimensionality of the
 algebras $\mathcal{B}_{\x}$ these compact spaces 
can be represented as $B_{\x}^{{\zz}^{\nu}}$
with finite sets $B_{\x}$ for all $\x\in A_{l,\eta}^{c}$.
By the Shannon-McMillan-Breiman theorem (cf. \cite{ornstein},
\cite{lindenstrauss})
\begin{equation}\label{eq:smb}
  \lim_{\Lambda(\n)\nearrow {\nn}^{\nu}}-\frac{1}{\#\Lambda(\n)}\log m_{\x}^{(\n)}(\omega_{\n})=h_{\x}
\end{equation}
$m_{\x}$-almost surely and in $L^{1}(m_{\x})$ for all $\x\in
A_{l,\eta}^{c}$, where $h_{\x}$ denotes the Kolmogorov-Sinai entropy
of $m_{\x}$, and $\omega_{\n} \in {B_{\x}}^{\Lambda(\n)}$ are the
  components of $\omega\in B_{\x}^{{\zz}^{\nu}}$ corresponding to 
the box $\Lambda (\n)$. Actually, as we shall see, we need the theorem
cited above only in its weaker form (convergence in probability) known as Shannon-McMillan theorem. 
For each $n$ and $\x\in A_{l,\eta}^{c}$ let
\begin{eqnarray}
 C_{\x}^{(\n)}& := & \{ \omega_{\n} \in B_{\x}^{ (\n)}|\ |-\frac{1}{\#\Lambda(\n)}\log
m_{\x}^{(\n)}(\omega_{\n}) -h_{\x}|<\delta\}\nonumber \\
& = & \{ \omega_{\n} \in B_{\x}^{ (\n)}|\ e^{-\#\Lambda(\n) \cdot (h_{\x}+\delta)}<m_{\x}^{(\n)}(\omega_{\n})
<e^{-\#\Lambda(\n) \cdot (h_{\x}-\delta)}\}.\nonumber 
\end{eqnarray}
Since lower bounds on the probability imply upper bounds on the cardinality
we obtain
\begin{equation}  \label{eq:trace}
 \# C_{\mathbf{x}}^{\left( \mathbf{n}\right) }=\textrm{tr}_{\mathbf{n}%
}\left( p_{\mathbf{x}}^{\left( \mathbf{n} \right)}\right) \leq e^{\#
\Lambda(\n) \cdot ( h_{\mathbf{x}}+\delta ) }\leq
e^{\# \Lambda( \mathbf{n}) \cdot \left( l^{\nu }\left( s\left(
\Psi \right) +\eta \right) +\delta \right) }
\end{equation}
where $ p_{\mathbf{x}}^{\left( \mathbf{n}\right) }$ is the projector
in ${\mathcal B}_{\x}^{(\n)} $ corresponding to the function
$1_{C_{\x}^{(\n)}} $. In the last inequality we have used that $h_{\x}\leq
S(\Psi_{\x}^{(\Lambda(l))}\upharpoonright \mathcal{B}_{\x})<l^{\nu}(s(\Psi)+\eta)$ for all $\x\in
A_{l,\eta}^{c}$ by
\[h_{\x}=\lim_{\Lambda(\n)\nearrow {\nn}^{\nu}}\frac{1}{\#\Lambda(\n)}H(m_{\x}^{(\n)})=
\inf_{\Lambda(\n)}\frac{1}{\#\Lambda(\n)}H(m_{\x}^{(\n)}),\]
 (cf. ~\cite{ruelle}), and by
(\ref{eq:goodentr}). Here $H$ denotes the Shannon entropy.\\
From (\ref{eq:smb}) it follows that there is an  $N\in{\nn}$
(depending on $l$) such that
for all $\n\in {\nn}^{\nu}$ with $\Lambda(\n)\supset \Lambda(N)$
\begin{equation} \label{eq:typ}
  m_{\x}^{(\n)}(C_{\x}^{(\n)})\geq 1-\frac{\varepsilon}{2},\qquad \forall
 \ \x\in A_{l,\eta}^{c}.
\end{equation}
For each $\y \in{\nn}^{\nu}$ with $y_{i}\geq Nl$ let
$y_{i}=n_{i}l+j_{i}$, 
where $n_{i}\geq N$ and $0\leq j_{i}<l$. We set 
\[
q:=\bigvee_{\x\in A_{l,\eta}^{c}}p_{\x}^{(\n)}. 
\] 
and denote by $q_{\y}$ the embedding of $q$ in $ \mathcal{A}^{(\y)}$. By
 (\ref{eq:typ}) and (\ref{eq:typical}) we obtain
\begin{eqnarray}
\Psi^{(\y)}(q_{\y})&=&\frac{1}{\#\Lambda(\mathbf{k}(l))}\sum_{\x\in
  \Lambda(\mathbf{k}(l))} \Psi_{\x}^{(\y)}(q_{\y})
 \nonumber \\
&\geq&\frac{1}{\#\Lambda(\mathbf{k}(l))}\#
A_{l,\eta}^{c} \cdot (1-\frac{\varepsilon}{2})\geq (1-\varepsilon).\nonumber
\end{eqnarray}
Thus the condition in the definition of $\beta_{\varepsilon,\y}(\Psi)$ is satisfied.
Moreover by (\ref{eq:trace})
\begin{eqnarray}
 \beta_{\varepsilon,\y}(\Psi) & \leq &
 \log
\textrm{tr}_{\y}(q_{\y})\nonumber \\
& \leq  & \log\sum_{\x\in
  A_{l,\eta}^{c}}e^{\#\Lambda(\n) \cdot (h_{\x}+\delta)}+ \textrm{tr}_{\Lambda(\y)\setminus\Lambda(l\n)}\idn  \nonumber \\
& \leq  & \log\#A_{l,\eta}^{c} \cdot 
e^{\#\Lambda(\n)(l^{\nu}(s(\Psi)+\eta)+\delta)})
 + \textrm{tr}_{\Lambda(\y)\setminus\Lambda(l\n)}\idn \nonumber \\
& \leq  & \log\#A_{l,\eta}^{c} +
\#\Lambda(\n) (l^{\nu}(s(\Psi)+\eta)+
\delta) \nonumber \\
& & + \textrm{tr}_{\Lambda(\y)\setminus\Lambda(l\n)}\idn \nonumber \\
& \leq &
\log\#A_{l,\eta}^{c} + \#(\Lambda(\y))(s(\Psi)+\eta+\delta)\nonumber
\\
& & + \textrm{tr}_{\Lambda(\y)\setminus\Lambda(l\n)}\idn. \nonumber 
\end{eqnarray}
We can conclude from this that 
\[ \limsup_{\Lambda(\y) \nearrow \nn^{\nu}}\frac{1}{\#(\Lambda(\y)) }\beta_{\varepsilon,\y}(\Psi)\leq s(\Psi)+\eta+\delta,\] 
because $\#A_{l,\eta}^{c}$ does not depend on $\n$ and
$\Lambda(\y)\nearrow {\nn}^{\nu}$ 
 if and only if $\Lambda(\n)\nearrow {\nn}^{\nu}$. This leads to
\[\limsup_{\Lambda(\y)\nearrow {\nn}^{\nu} }\frac{1}{\#\Lambda(\y)}\beta_{\varepsilon,\y}(\Psi)\leq s(\Psi),\]
since $\eta$, $\delta>0$ were chosen arbitrarily. $\qquad \qed $\\
%

Let $\nu \in \nn$. For $\n  = (n_{1}, \dots, n_{\nu}) \in \nn^{\nu}$
we define $|\n| := \prod_{i=1}^{\nu}n_{i}$ and write
 $\n \rightarrow \infty$ alternatively for  $ \Lambda(\n) \nearrow
 \nn^{\nu}$. Further we
 introduce the notation 
\begin{eqnarray*}
\n \geq \mathbf{m}\ :\Longleftrightarrow \ n_{i} \geq m_{i}, \quad
\forall i \in \{1, \dots, \nu\}  .
\end{eqnarray*}
Recall that for a probability
distribution $P$ on a finite set $A$ the Shannon entropy is defined by 
\begin{eqnarray*}
H(P):= -\sum_{a \in A} P(a) \log
P(a).
\end{eqnarray*}
\begin{lemma}\label{class_resultat}
Let  $D>0$ and $ \{(A^{(\n)}, P^{(\n)})\}_{\n \in \nn^{\nu}} $ be a family,
where each  $ A^{(\n)}$ is a finite set with $\frac{1}{|\n|}\log
\#A^{(\n)}\leq D$ for all ${\mathbf{n}} \in \nn^{\nu}$ and
$P^{(\n)}$ is a probability
distribution on $A^{(\n)}$. Define
\begin{eqnarray}\label{alpha}
\alpha_{\eps,\n}(P^{(\n)}):= \min \{ \log  \# \Omega |\ \Omega \subset A^{(\n)}, P^{(\n)}(\Omega) \geq  1-
  \eps \}.
\end{eqnarray}
If $ \{( A^{(\n)}, P^{(\n)})\}_{\n \in \nn^{\nu}} $ satisfies the following two conditions:

\begin{enumerate}

\item\label{Ann1}
$\lim_{\n \to \infty} \frac{1}{|\n|} H(P^{(\n)}) = h < \infty$

\item\label{Ann2} 
$\limsup_{\n \to \infty} \frac{1}{| \n |}\alpha_{\eps,\n}(P^{(\n)}) \leq h, \qquad \forall \eps \in
(0,1)$

\end{enumerate}
then for every $\eps \in (0,1)$
\begin{eqnarray}\label{alpha=h}
\lim_{\n \to \infty} \frac{1}{|\n|} \alpha_{\eps, \n}(P^{(\n)})=h.
\end{eqnarray}
\end{lemma}
Note that we do not expect either $\{A^{(\n)}\}_{\n \in \nn^{\nu}}$ or $\{P^{(\n)}\}_{\n \in \nn^{\nu}}$ to fulfill any
consistency conditions. We will see later on that this is the
important point for why Lemma \ref{class_resultat} will be useful in the
non-commutative setting.\\
\\
\textbf{Proof of Lemma \ref{class_resultat}:}
Let $\delta >0$ and define
\begin{eqnarray}
A_{1}^{\left( \mathbf{n}\right) }(\delta)&:=& \left\{ a\in A^{\left( \mathbf{n}\right)
}\mid P^{\left( \mathbf{n}\right) }\left( a\right) >e^{-\left| \mathbf{n}%
\right| \left( h-\delta \right) }\right\}, 
\nonumber \\
A_{2}^{( \mathbf{n}) }(\delta)&:=& \left\{ a\in A^{\left( \mathbf{n}\right)
}\mid e^{-\left| \mathbf{n}\right| \left( h+\delta \right) }\leq P^{\left( 
\mathbf{n}\right) }\left( a\right) \leq e^{-\left| \mathbf{n}\right| \left(
h-\delta \right) }\right\}, \nonumber \\
A_{3}^{\left( \mathbf{n}\right) }(\delta)&:=& \left\{ a\in A^{\left( \mathbf{n}\right)
}\mid P^{\left( \mathbf{n}\right) }\left( a\right) <e^{-\left| \mathbf{n}%
\right| \left( h+\delta \right) }\right\}. 
\nonumber
\end{eqnarray}
We fix an arbitrary $\delta>0$ and use the abbreviation
$A^{(\n)}_{i}=A^{(\n)}_{i}(\delta)$, $i \in \{1,2,3\}$. To see that $\lim\limits_{\mathbf{n}\rightarrow \infty }P^{\left( \mathbf{n}%
\right) }\left( A_{3}^{\left( \mathbf{n}\right) }\right) =0$ assume the
contrary and observe that the upper bound on the probability of elements
from $A_{3}^{\left( \mathbf{n}\right) }$ implies a lower bound on the
cardinality of elements from $A_{3}^{\left( \mathbf{n}\right) }$ needed to
cover say an $\varepsilon-$fraction, $\eps \in (0,1)$, of $A^{\left( \mathbf{n}\right) }$
with respect to $P^{\left( \mathbf{n}\right) }$. Namely one has $\min \left\{ \#
C\mid C\subset A_{3}^{\left( \mathbf{n}\right) },P^{\left( \mathbf{n}\right)
}\left( C\right) >\varepsilon \right\}
>\varepsilon \cdot e^{\left| \mathbf{n}\right| \left( h+\delta \right) }$ which
contradicts condition 2 in the lemma. Furthermore the set $A_{3}^{\left( 
\mathbf{n}\right) }$ cannot asymptotically contribute to the
mean entropy $h$
since
\begin{eqnarray*}
-\frac{1}{\left| \mathbf{n}\right| }\sum\limits_{a\in A_{3}^{\left( \mathbf{%
n}\right) }}P^{\left( \mathbf{n}\right) }\left( a\right) \log P^{\left( 
\mathbf{n}\right) }\left( a\right)\\ \leq  - \frac{1}{\left| \mathbf{n}\right| }\sum\limits_{a\in A_{3}^{\left(
\mathbf{n}\right) }}P^{\left( \mathbf{n}\right) }\left( a\right) \log \frac{1%
}{\# A_{3}^{\left( \mathbf{n}\right) }}P^{\n}(A_{3}^{\n})
\end{eqnarray*}
and 
\begin{eqnarray*}
& & \lim_{\n \rightarrow \infty } -\frac{1}{%
\left| \n \right| } \sum\limits_{a\in A_{3}^{\left( \mathbf{n}\right)
}}P^{\left( \mathbf{n}\right) }\left( a\right) \log \frac{1}{\#
A_{3}^{\left( \mathbf{n}\right) }}P^{(\n)}(A_{3}^{(\n)}) \\ &=& \lim_{\n
\to \infty}
\frac{1}{|\n|}\left(P^{(\n)}(A^{(\n)}_{3}) \log \# A^{(\n)}_{3}- P^{(\n)}(A^{(\n)}_{3})\log
P^{(\n)}(A^{(\n)}_{3}) \right)=0.
\end{eqnarray*}
Here we used the
fact that $\frac{\log \# A^{\left( \mathbf{n}\right) }}{\left| \mathbf{n}%
\right| }$ stays bounded and $-\sum p_{i}\log p_{i}\leq -\sum p_{i}\log q_{i}
$ for finite vectors $\left( p_{i}\right) ,\left(
q_{i}\right)$ with $\sum_{i}p_{i}= \sum_{i}q_{i} \leq 1$ and
$p_{i},q_{i}\geq 0$. Since $A_{3}^{\left( \mathbf{n}\right) }$ does not
contribute to the entropy one easily concludes that $\lim\limits_{%
\mathbf{n}\rightarrow \infty }P^{\left( \mathbf{n}\right) }\left(
A_{1}^{\left( \mathbf{n}\right) }\right) =0$ because otherwise $%
\liminf_{ \n \rightarrow \infty } \frac{1}{%
\left| \n \right| }H\left( P^{\left( \mathbf{n}\right) }\right) <h$ would
hold. Recall that $\delta>0$ was chosen arbitrarily. \\
Thus
\begin{eqnarray}\label{a2-probab}
\lim\limits_{\mathbf{n}\rightarrow \infty }P^{\left( \mathbf{n}\right)
}\left( A_{2}^{\left( \mathbf{n}\right) }(\delta)\right) =1, \qquad
\forall \delta>0.
\end{eqnarray}
Consequently the lemma follows
since $P^{\left( \mathbf{n}\right) }\left( \Omega \right) \geq 1-\varepsilon 
$ implies $P^{\left( \mathbf{n}\right) }\left( \Omega \cap A_{2}^{\left( 
\mathbf{n}\right) }(\delta)\right) \geq (1-\varepsilon)^{2} $ for $\left| \mathbf{n}%
\right| $ sufficiently large and one needs \ at least $\left( 1-\varepsilon
\right)^{2}\cdot e^{\left| \mathbf{n}\right| \left( h-\delta \right) }$ elements from 
$A_{2}^{\left( \mathbf{n}\right) }(\delta)$ to cover $\Omega \cap A_{2}^{\left( 
\mathbf{n}\right) }(\delta)$ and $\delta $ can be chosen arbitrarily
small. $\qquad \qed$
\\ \\
%
%
%
%
%
%
\textbf{Proof of Proposition \ref{beta=s}:}
${\cal A}^{(\n)}$ as a finite dimensional $C^{*}$-algebra is isomorphic to a
finite direct sum $\bigoplus_{j=1}^{M}{\cal B}(\hr_{j}^{(\n)})$, where each
$\hr_{j}^{(\n)}$ is a Hilbert space with $\dim \hr_{j}^{(\n)}= d_{j}^{(\n)}< \infty$ and any minimal projector in ${\cal
  A}^{(\n)}$ is represented by a one-dimensional projector on
$\hr^{(\n)}:=\bigoplus_{j=1}^{M} \hr_{j}^{(\n)}$ with $\dim \hr^{(\n)}=\sum_{j=1}^{M}d_{j}^{(\n)}=:d_{\n}$. Note that $\bigoplus_{j=1}^{M}{\cal
  B}(\hr_{j}^{(\n)}) \subset {\cal B}(\hr^{(\n)})$.
Consider the spectral representation of the density operator $D_{\n}$ of
$\Psi^{(\n)}$ in ${\cal B}(\hr^{(\n)})$:
\begin{eqnarray*}
D_{\n}= \sum_{i=1}^{d_{\n}}
\lambda_{i}^{(\n)} |q^{(\n)}_{i} \rangle \langle q^{(\n)}_{i} |, \qquad \lambda_{i}^{(\n)} \in [0,1],\
| q^{(\n)}_{i} \rangle \in \hr^{(\n)}.
\end{eqnarray*}
For $\n=(n_{1}, \dots, n_{\nu}) \in \nn^{\nu}$ let $A^{(\n)}$ be the finite set consisting
of the eigen-projectors $q^{(\n)}_{i} := |q^{(\n)}_{i} \rangle \langle
q^{(\n)}_{i} | $ of $\Psi^{(\n)}$, i.e.
\begin{eqnarray}\label{alphabet}
A^{(\n)}:=
\{q^{(\n)}_{i} \}_{i=1}^{d_{\n}}.
\end{eqnarray}
Let $P^{(\n)}$ be the probability
distribution on $A^{(\n)}$ given by:
\begin{eqnarray}\label{probability}
P^{(\n)}(q^{(\n)}_{i}):=\Psi^{(\n)}(q^{(\n)}_{i})=\lambda^{(\n)}_{i}.
\end{eqnarray}
Recall that $|\n|=\prod_{i=1}^{\nu}n_{i}$. Let $D:= \log (\dim
\hr^{(0)})$, then $\frac{1}{|\n|}\log \# A^{(\n)} \leq D$ for all $\n \in \nn^{\nu}$. We show that the family $\{(A^{(\n)}, P^{(\n)})\}_{\n \in \nn^{\nu}}$ fulfills both
conditions in Lemma \ref{class_resultat} and consequently
\begin{eqnarray}\label{alpha=hicks}
\lim_{ \n \to
  \infty}\frac{1}{|\n|}\alpha_{\eps,\n}(P^{(\n)})=\lim_{\n \to \infty} \frac{1}{|\n|}H(P^{(\n)}), \qquad \forall \eps
\in (0,1).
\end{eqnarray}
It is clear that $ H(P^{(\n)})= -\sum_{i=1}^{d_{\n}} \lambda^{(\n)}_{i}
\log \lambda^{(\n)}_{i} =S(\Psi^{(\n)}) $. Thus
\begin{eqnarray}\label{h=s}
h:=\lim_{\n \to \infty} \frac{1}{|\n|}H(P^{(\n)})=s(\Psi).
\end{eqnarray}
Next assume the following ordering:
\begin{eqnarray}
i<j \ \Longrightarrow \ \lambda^{(\n)}_{i} \geq \lambda^{(\n)}_{j} \nonumber
\end{eqnarray}
and define for $\eps \in (0,1)$
\begin{eqnarray}
n_{\eps, \n}:= \min \{k \in \{1, \dots, d_{\n}\} | \ \sum^{k}_{j=1} \nonumber
\lambda^{(\n)}_{j}\geq 1-\eps \}.
\end{eqnarray}
Thus $ \alpha_{\eps, \n}(P^{(\n)})= \log \#(
\{q^{(\n)}_{i}\}_{i=1}^{n_{\eps, \n}})= \log n_{\eps, \n}$. We claim :
\begin{eqnarray}\label{beta=alpha}
\alpha_{\eps, \n}(P^{(\n)})=\beta_{\eps, \n}(\Psi^{(\n)}) , \qquad \forall
\eps \in (0,1).
\end{eqnarray}
From $\Psi^{(\n)}(\sum_{i =1}^{n_{\eps,\n}}q^{(\n)}_{i}) \geq
1-\eps$ and $\textrm{tr}_{\n}\sum_{i =1}^{n_{\eps,
    \n}}q^{(\n)}_{i}=n_{\eps, \n} $ it is obvious that  $\beta_{\eps, \n}(\Psi^{(\n)}) \leq \alpha_{\eps, \n}(P^{(\n)})$.\\
Assume $\beta_{\eps, \n}(\Psi^{(\n)}) < \alpha_{\eps, \n}(P^{(\n)}) $. Then
there exists a projector $q \in {\cal A}^{(\n)}$ with $\Psi^{(\n)}(q)\geq 1-\eps$ such that
$m:=\textrm{tr}_{\n}\ q < n_{\eps, \n}$. Let $\sum_{i=1}^{m} |
q_{i} \rangle \langle q_{i} |$, where $|
q_{i} \rangle \in \hr^{(\n)}$, be the spectral representation
of $q$. For $D_{\n}$ as density matrix on $\hr^{(\n)}$ we use Ky
Fan's maximum principle, \cite{bhatia}, and obtain the contradiction
\begin{eqnarray}
1- \eps \leq \Psi^{(\n)}(q)=\textrm{tr}_{\n}D_{\n}q =
\sum_{i=1}^{m}\langle q_{i}, D_{\n}q_{i} \rangle  \leq
\sum_{i=1}^{m} \lambda^{(\n)}_{i}< 1- \eps. \nonumber
\end{eqnarray}
$\Psi$ is ergodic. Thus we can apply Lemma \ref{lim_beta<s}:
\begin{eqnarray}\label{lim_sup}
\limsup_{\Lambda(\n) \nearrow \nn^{\nu}}\frac{1}{\#(\Lambda(\n))}\beta_{\eps,\n}(\Psi^{(\n)})\leq s(\Psi), \qquad
\forall \eps \in (0,1).
\end{eqnarray}
Setting (\ref{beta=alpha}) and (\ref{h=s}) in (\ref{lim_sup}) and
using that $\#(\Lambda(\n))=|\n|$  we obtain
\begin{eqnarray}\label{2nd_assump}
\limsup_{\n \to \infty} \frac{1}{|\n|}\alpha_{\eps,\n}(P^{(\n)})\leq h, \qquad \forall
\eps \in (0,1).
\end{eqnarray}
With (\ref{h=s}) and (\ref{2nd_assump}) both conditions in
Lemma \ref{class_resultat} are satisfied. It follows (\ref{alpha=hicks}). Now we set back (\ref{beta=alpha})
and (\ref{h=s}) in (\ref{alpha=hicks}) and arrive at \[ \lim_{\Lambda(\n)
  \nearrow \nn^{\nu}}\frac{1}{\#(\Lambda(\n))}\beta_{\eps,\n}(\Psi)=s(\Psi), \qquad \forall \eps
\in (0,1). \qquad \qed\]
%
%
%
%
%
\textbf{Proof of the Quantum Shannon-McMillan Theorem:}\\
Fix $\delta >0.$ Adopt the family $\{(A^{(\mathbf{n})},P^{(%
\mathbf{n)}})\}_{\mathbf{n\in }\mathbb{N}^{\nu }}$ and further
notations from the proof of Proposition  \ref{beta=s}. Choose some $\delta ^{\prime }<\delta .$ Let $A_{2}^{(\mathbf{n)}%
}(\delta ^{\prime })$ be the subset of $A^{(\mathbf{n})}$ defined in the
proof of Lemma \ref{class_resultat} with $h=s(\Psi )$, appropriate to
(\ref{h=s}). Let $I_{\n}(\delta^{\prime }):=\{ i \in \{1, \dots, d_{\n}\}  |\ q_{i}^{\n}\in A_{2}^{(\n)
}(\delta ^{\prime }) \}$. Set
\begin{eqnarray*}
p_{\mathbf{n}}(\delta)=\sum\limits_{I_{\n}(\delta ^{\prime })}q_{i}^{\mathbf{n}}.
\end{eqnarray*}
By (\ref{a2-probab}) there exists an $\mathbf{N}_{\delta} \in \nn^{\nu}$ such
that $p_{\mathbf{n}}(\delta)$ is a projector with
\begin{eqnarray*}
\Psi^{(\n)}(p_{\mathbf{n}}(\delta))=P^{(\mathbf{n})}(A_{2}^{(\mathbf{n)}}(\delta ^{\prime }))\geq
1-\delta, \qquad \forall \n \geq \mathbf{N}_{\delta}.
\end{eqnarray*}
Any minimal projector $0 \not= p \in {\cal A}^{(\n)}$ with $p
\leq p_{\mathbf{n}%
}(\delta)$ is represented as a one-dimensional projector $| p
\rangle \langle p |$ on $\hr^{(\n)}$, such that $| p
\rangle= \sum_{i \in I_{\n}(\delta^{\prime})}\gamma
_{i}|q_{i}^{\mathbf{n}} \rangle $ and $\sum_{i \in I_{\n}(\delta^{\prime})}\gamma^{2}_{i}=1$. Hence
\begin{eqnarray*}
\Psi ^{(\mathbf{n})}(p)=\sum\limits_{i \in I_{\n}(\delta ^{\prime })}\gamma _{i}^{2}\lambda _{i}^{(\mathbf{n})}
\end{eqnarray*}
is a weighted average of the eigenvalues $\lambda _{i}^{(\mathbf{n})}$ corresponding to the set $A_{2}^{(\mathbf{n)%
}}(\delta ^{\prime })$. Thus we obtain by the definition of this set
\begin{eqnarray}\label{item2}
e^{-\#\Lambda (\mathbf{n})(s(\Psi )+\delta )}<\Psi ^{(\mathbf{n}%
)}(p)<e^{-\#\Lambda (\mathbf{n})(s(\Psi )-\delta )}.
\end{eqnarray}
Using the
linearity of $\Psi^{(\n)}$ and applying (\ref{item2}) to the
projectors $q^{(\n)}_{i}$ we arrive at
the following estimation
\begin{eqnarray*}
e^{\# \Lambda(\n)(s(\Psi)-\delta)}< \textrm{tr}_{\n} p_{\n}(\delta)<
e^{\#   \Lambda (\n)(s(\Psi) + \delta)},
\end{eqnarray*}
if $\n$ is large enough.
We have shown all assertions of the theorem. $\qquad \qed$
\section{Comment}  \label{conclusions}
It seems that our result can be extended to the case of discrete
amenable group actions on quasilocal $C^{\ast }-$algebras. Limits in
this situation have to be taken along tempered F\o lner sequences. The
relevant classical theorem in this setting can be found in \cite{lindenstrauss}.

\begin{acknowledgement}
We want to express our deep thanks to Ruedi Seiler who
emphasized the relevance of a quantum version of the
Shannon-McMillan theorem to us. He stimulated our work on the problem in
many discussions and with helpful comments. Without his continuous
encouragement and optimism this paper would never have been
written. \\ \\
This work was supported by the DFG via the SFB 288 ``Quantenphysik und
Differentialgeometrie'' at the TU Berlin and the Forschergruppe
``Stochastische Analysis und gro\ss e Abweichungen'' at University of Bielefeld.
\end{acknowledgement}


\begin{thebibliography}{99}
\bibitem{benatti}
F. Benatti, T. Hudetz, A. Knauf, Quantum Chaos and Dynamical Entropy,
Commun. Math. Phys. 198, 607-688 (1998)

\bibitem{bhatia}
R. Bhatia, Matrix Analysis, Graduate Texts in Mathematics 169,
Springer, New York 1997

\bibitem{bksss2}
I. Bjelakovi\'c, T. Kr\"uger, Ra. Siegmund-Schultze, A. Szko\l a,
 Chained Typical Subspaces - a Quantum Version of Breiman's Theorem, in preparation

\bibitem{bratteli}
O.Bratteli, D.W. Robinson, Operator Algebras and Quantum Statistical
Mechanics I, Springer, New York 1979

\bibitem{breiman}
L. Breiman, The Individual Ergodic Theorem of Information Theory,
Ann. Math. Stat. 28, Nr. 3, 809-811 (1957) 

\bibitem{datta}
N. Datta, Y. Suhov, Data Compression Limit for an Information Source
of Interacting Qubits, math-phys/0207069

\bibitem{petz}
F. Hiai, D. Petz, The Proper Formula for Relative
  Entropy and its Asymptotics in Quantum Probability,
  Commun. Math. Phys. 143, 99-114 (1991)

\bibitem{hiai}
F. Hiai, D. Petz, Entropy Densities for Algebraic States, J. 
Functional Anal. 125, 287-308 (1994)

\bibitem{gibbs}
F. Hiai, D. Petz, Entropy Densities for Gibbs States of Quantum Spin
Systems, Rev. Math. Phys. 5 No.4, 693-712 (1993)

\bibitem{israel}
R.B. Israel, Convexity in the Theory of Lattice Gases,
Princeton, New Jersey 1979

\bibitem{suhov}
O. Johnson, Y. Suhov, The von Neumann Entropy and Information Rate
for Ideal Quantum Gibbs Ensembles, math-phys/0109023

\bibitem{jozsa}
R.Jozsa, B. Schumacher, A New Proof of the Quantum Noiseless Coding
Theorem, Journal of Modern Optics Vol.41, No.12, 2343-2349 (1994)

\bibitem{king} C. King, A. Le\'sniewski, Quantum Sources and a
  Quantum Coding Theorem, J. Math. Phys. 39 (1), 88-101 (1998) 

\bibitem{lindenstrauss}
E. Lindenstrauss, Pointwise Theorems for Amenable Groups, 
Invent. Math. 146, 259-295 (2001)

\bibitem{neshveyev}
S.Neshveyev, E.St\o rmer, The McMillan Theorem for a Class of
Asymptotically Abelian $C^{*}$-Algebras, Ergod. Th. \& Dynam.
Sys. 22, 889-897 (2002) 

\bibitem{niel} 
M. A. Nielsen, I. L. Chuang, Quantum Computation and
 Quantum Information, Cambridge University Press, Cambridge 2000

\bibitem{ohya}
M. Ohya, D. Petz, Quantum Entropy and its Use, Springer, Berlin 1993

\bibitem{ornstein}
D. Ornstein, B. Weiss, The Shannon-McMillan-Breiman Theorem
for a Class of Amenable Groups, Israel J. Math. 44 (1), 53-60
 (1983)

\bibitem{mosonyi}
D. Petz, M. Mosonyi, Stationary Quantum Source Coding,
J. Math. Phys. 42, 4857-4864 (2001)

\bibitem{shields}
P. C. Shields, The Ergodic Theory of Discrete Sample Paths, Graduate
Studies in Mathematics Vol. 13, American Mathematical Society 1996    

\bibitem{ruelle} 
D. Ruelle, Statistical Mechanics, W.A. Benjamin, New York 1969





\end{thebibliography}
\end{document}